\documentclass{article}

\usepackage{amssymb}
\usepackage{amsmath}

\newtheorem{theorem}{Theorem}
\newtheorem{lemma}{Lemma}[section]
\newtheorem{propo}{Proposition}[section]

\newtheorem{defin}{Definition}[section]
\newtheorem{remark}{Remark}[section]

\newcommand{\ep}{\varepsilon}

\newcommand{\de}{\delta}
\newcommand{\al}{\alpha}

\newcommand{\e}{\ep}
\newcommand{\ad}{d^2+1 \choose 2}

\newcommand{\N}{\mathbb{N}}
\newcommand{\R}{\mathbb{R}}

\newcommand{\Grd}[1][d]{\textrm{\textbf{Gr}}_{#1}}
\newcommand{\Gd}[1][d]{\textrm{\textbf{Graph}}_{#1}}
\newcommand{\CGd}[1][d]{\textrm{\textbf{C-Graph}}_{#1}}
\newcommand{\CGrd}[1][d]{\textrm{\textbf{CGr}}_{#1}}
\newcommand{\BGrd}[1][d]{\textrm{\textbf{BGr}}_{#1}}
\newcommand{\hb}{\hat{\BGrd}}

\newcommand{\cgim}{CG_n^{i,M}}

\newcommand{\defeq}{\stackrel{\mbox{\scriptsize def}}{=}}

\newcommand{\diam}{\textrm{diam}}

\newcommand{\G}{\textrm{\textbf{G}}}

\newcommand{\CG}[1][G]{\textrm{C\textbf{#1}}}

\newcommand{\I}{\mathcal{I}}
\newcommand{\EI}{\mathcal{EI}}
\newcommand{\Gseq}[1][G]{\{ #1_n \}_{n=1}^{\infty}}

\def\proof{\smallskip\noindent{\it Proof.} }
\newcommand{\qed} {\hspace {0.1in} \rule {1.5mm} {3.5mm}}
\title{An analogue of the Szemeredi Regularity Lemma for
bounded degree graphs\footnote{AMS
Subject Classification: Primary 05C99, Secondary 37A20,
\, Research sponsored by OTKA Grant No. 67867, No. 69062}}

\author{G\'abor Elek, G\'abor Lippner}

\begin{document}
\maketitle

\begin{abstract}
We show that a sufficiently large graph of bounded degree can be
decomposed into quasihomogeneous pieces. 
The result can be viewed as a ``finitarization''
of the classical Farrell-Varadarajan Ergodic Decomposition Theorem.
\end{abstract}

\section{Introduction}

 In order to state our result we need to recall some basic definitions.
Let $\Gd$ denote the set of all connected finite simple
 graphs $G$ (up to isomorphism)
  for which $\deg(x) \leq d$ for every $x \in V(G)$.
For a graph $G$ and $x,y \in V(G)$ let $d_G(x,y)$ denote the distance
of $x$ and $y$, that is the length of the shortest path from $x$ to
$y$.  A rooted $(r,d)$-ball is a graph $G \in \Gd$ with a marked
vertex $x \in V(G)$ called the root such that $d_G(x,y) \leq r$ for
every $y \in V(G)$.  By $U^{r,d}$ we shall denote the set of rooted
$(r,d)$-balls (up to rooted isomorphism).

If $G \in \Gd$ is a graph and $x\in V(G)$ then $B_r(x) \in U^{r,d}$
shall denote the rooted $(r,d)$-ball around $x$ in $G$.
For any $\alpha \in U^{r,d}$ and $G \in \Gd$ we define the set
  $T(G,\alpha) \defeq \{x \in V(G): B_r(x) \cong \alpha\}$ and let
  $p_G(\alpha) \defeq \frac{|T(G,\alpha)|}{|V(G)|}$.
Fix an enumeration of all possible $(r,d)$-balls, for every $r\geq 1$:
$(\alpha_1,\alpha_2,\dots)$. Let us define the {\bf statistical distance}
of two graphs $G$ and $H$ as
$$d_s(G,H)=\sum^\infty_{i=1}\frac{1}{2^i} |p_G(\alpha_i)-p_H(\alpha_i)|\,.$$
It is easy to see that $d_s(G,H)$ defines a metric on $\Gd$.
We define $d_s$ for not necessary connected graphs as well. In this
case $d_s$ defines a pseudo-distance.
\begin{remark} \label{rem1}
Note that if $F_1$, $F_2$,\dots $F_k$ are finite connected graphs
then for any $\epsilon>0$ there exists a $\delta>0$ such that if
$d_s(G,H)\leq\delta$ then for any $1\leq i \leq k$ 
$|\mbox{dens}(F_i,G)-\mbox{dens}(F_i,H)|\leq\epsilon$, where
$$\mbox{dens}(F_i,G)=\frac{\mbox{the number of subgraphs of $G$
isomorphic to $F_i$}}{|V(G)|}\,.$$
The number $\mbox{dens}(F_i,G)$ is called the ``sparse'' $F_i$-density of $G$.
\end{remark}

\noindent
If $J\subset G$ is a spanned subgraph then $E(J, G\backslash J)$ denotes
the number of edges between the vertices of $J$ and $G\backslash J$.
The following is our key definition.
\begin{defin}
$G\in\Gd$ is called $(\epsilon,\lambda,\delta)$-quasihomogeneous, if
for any spanned subgraph $H\subset G$ such that
\begin{itemize}
\item $\lambda\leq \frac{|V(H)|}{|V(G)|}$
\item $E(H, G\backslash H)\leq \epsilon|V(G)|$
\end{itemize}
we have $d_s(G,H)\leq \delta$.
\end{defin}
Informally speaking $G$ is "quasihomogeneous" if for any large enough
spanned graph $H$ which is connected to $G\backslash H$ by only a small
amount of edges, the subgraph densities of $G$ and $H$ are 
very similar.
Now let us recall the regularity lemma of Szemeredi for
dense graphs. Let $G$ be a graph and
$X\subset V(G)$, $Y\subset V(G)$ be disjoint subsets. The density
of the pair $X$,$Y$ is
$$p(X,Y):=\frac{|E(X,Y)|}{|X||Y|}\,.$$
The pair $X,Y$ is called $\epsilon$-quasirandom if for any subsets
$A\subset X$, $B\subset Y$, $|A|\geq\epsilon |X|$, $|B|\geq\epsilon |Y|$
$$|p(A,B)-p(X,Y)|<\epsilon\,.$$
Now we have a similar subgraph counting principle as in Remark \ref{rem1}.
If $F_1$, $F_2$,\dots, $F_k$ are fixed simple connected graphs
then for any $\gamma>0$ there exists a $\epsilon>0$ such that if
$X,Y$ is a $\epsilon$-quasirandom pair, and $A$,$B$ are subsets
as above then for any $1\leq i \leq k$:
$$|\mbox{density}(F_i,E(X,Y))-\mbox{density}(F_i,E(A,B))|\leq\gamma\,,$$
where
$$\mbox{density}(F_i,G)=\frac{\mbox{the number of subgraphs of $G$
isomorphic to $F_i$}}{|V(G)|^{|F_i|}}\,.$$
The number $\mbox{density}(F_i,G)$ is called the $F_i$-density of $G$.
That is quasirandomness also implies a certain kind of quasihomogeneity.
According to the Szemeredi Regularity Lemma for any $\epsilon>0$
there exists $K(\epsilon)>\frac{1}{\e}$ and $N(\epsilon)>1$ such that for any 
graph $G$ with $|V(G)|\geq N(\e)$ one can remove $\epsilon |V(G)|^2$ edges
from $G$ such that the vertices of the
 remaining graph $G'$ can be partitioned into
$\frac{1}{\e}\leq K \leq K(\e)$ parts (of almost equal size)
 $V(G')=V(G'_1)\cup V(G'_2)\cup\dots\cup V(G'_K)$ and all pairs 
$(G'_i,G'_j)$ are
$\epsilon$-quasirandom. Our main theorem might be considered as a bounded
degree analogue of the Szemeredi Regularity Lemma:
\begin{theorem}~\label{mainthm}
For every $\delta>0,\lambda>0$ 
there exist positive integers $K
= K(\de,\lambda)$, $N = N(\de,\lambda)$ and a positive
constant $\ep=\ep(\delta,\lambda)<\delta$
for 
which the following hold: 
\begin{enumerate}\item
If $G\in \Gd$ is a finite connected graph that
has at least $N$ vertices, then it can be partitioned into $K+1$
disjoint graphs $G=G_1\cup G_2\cup\dots G_K\cup G_{\emptyset}$
by deleting at most $\de |E(G)|$ edges such that
\begin{itemize}
\item $G_{\emptyset}$ is an edge-less graph, $\frac{|V(G_\emptyset)|}
{|V(G)|}<\de$
\item Either $\frac{|V(G_i)|}
{|V(G)|}>\frac{\delta^2}{10\,d\,K}$ or $V(G_i)$ is empty.
\item All non-empty parts are $(\epsilon,\lambda,\delta)$-quasihomogeneous.
\end{itemize}
\item For any $\sigma > 0$ there is a $\tau > 0$ such that whenever
  $d_s(G,H) < \tau$ then there exist partitions for which
  $d_s(G_i,H_i) < \sigma$ and $\left|\frac{|V(G_i)|}{|V(G)|} - \frac{|V(H_i)|}
{|V(H)|}\right| < \sigma$ for all $i$'s where either $G_i$ or $H_i$ is
  non-empty.
\end{enumerate}
\end{theorem}
The second statement of the theorem can be interpreted that if two graphs
are close in terms of ``sparse'' subgraph densities then they
have similar quasihomogeneous partitions. One should note that according to
the result of Borgs, Chayes, Lovasz, T.S\'os and Vesztergombi \cite{BCh2}
if two dense graphs are close in terms of their subgraphs densities
then they have similar Szemeredi partitions.
\noindent
The proof of the theorem is based on a ``finitarization'' of the
Farrell-Varadarajan Ergodic Decomposition Theorem. The necessary background
on ergodic theory and its connections to graph theory shall be surveyed
in Section \ref{notationsec}. 
The proof of Theorem \ref{mainthm} shall be given
in Section \ref{mainsec}.

\section{Ergodic theory}~\label{notationsec}
\subsection{Borel equivalence relations
 and invariant measures}~\label{ergodicsec}

 In this subsection we recall some basic notions from Chapter I. of \cite{KM} 
on countable
Borel equivalence relations.
 Let $X$ be a compact metric space. $E\subset X\times X$ is a 
countable Borel equivalence relation if all the equivalence classes are
countable and $E$ is a Borel subset of $X\times X$. Typical example of a Borel 
equivalence relation is the orbit equivalence relation of a Borel action of
a countable group. As a matter of fact according to the theorem of Feldman and
Moore any Borel equivalence relation can be described this way.
A Borel probability measure $\mu$ on $X$ is $E$-invariant
 if its invariant under
a countable group action that defines $E$. Note that in this case $\mu$ is 
invariant under {\it all} the group action that defines $E$. Equivalently,
$\mu$ is invariant if for any Borel isomorphism $f:X\to X$ $f_*(\mu)=\mu$,
that is if $A\subseteq X$ is a Borel-set then $\mu(f^{-1}(A))=\mu(A)$.
The space of invariant probability measures is denoted by $\I_E$.
A measurable set $A\subseteq X$ is called $E$-invariant if for 
any $x\in A$ and 
$y\in X, x \sim_E y$: $y\in A$. 
The invariant measure $\mu$ is called {\bf ergodic}
if the $\mu$-measure of any invariant set is either $0$ or $1$. 
The space of ergodic probability measures is denoted by
$\EI_E$.
Note that the set of probability measures on $X$, $P(X)$ is compact convex
set of the topological vectorspace of all measures in the weak-topology 
(Banach-Alaoglu
Theorem). 
The space $\I_E$ is a convex subset of $P(X)$ and
$\EI_E$ can be identified as the set of extremal points in $\I_E$.

\noindent
Our main tool will be the following well known result (see
e{.}g{.}~\cite{KM}):

\begin{propo}[Ergodic Decomposition -- Farrell, Varadarajan]~\label{decomp} Let
$E$ be a countable Borel equivalence relation on $X$. Then $\I_E,
\EI_E$ are Borel sets in the standard Borel space $P(X)$ of
probability measures on $X$.  Now suppose $\I_E \neq \emptyset$.
Then $\EI_E \neq \emptyset$, and there is a Borel surjection $\pi :
X \to \EI_E$ such that
\begin{enumerate}
\item $\pi$ is $E$-invariant,
\item if $X_e = \{x : \pi(x) = e\}$, then $e(X_e) = 1$ and $E|X_e$ has a
unique invariant measure, namely $e$, and
\item if $\mu \in \I_E$, then $\mu = \int \pi(x) d\mu(x)$.
\end{enumerate}
 Moreover, $\pi$ is uniquely determined in the sense that, if $\pi'$
is another such map, then $\{x : \pi(x) \neq \pi'(x)\}$ is null with
respect to all measures in $\I_E$.
\end{propo}
We need a simple observation on the space $P(X)$. The space $P(X)$
is metrizable in the weak-topology by
$$d_X(\mu,\nu):=\sum^\infty_{n=1}\frac{1}{2^n}|(\mu-\nu)(f_n)|\,,$$
where $\{f_n\}^\infty_{n=1}$ is a countable dense set in the unit
ball of $C(X)$ (the space of continuous functions). Here $\mu(f)$ denotes
$\int_Xf\,d\mu$.
\begin{lemma} \label{l21}

\noindent
\begin{itemize}\item
If $\mu_n\to \mu$, $\nu_n\to\nu$ in $P(X)$ and $\lambda_n\to\lambda$ in $\R$
then $\sum_{n=1}^\infty \lambda_n\mu_n+(1-\lambda_n)\nu_n\to
\lambda\mu+(1-\lambda)\nu$.
\item If $T\subset P(X)$ is an arbitrary subset then if $p$ is in
the convex hull of $T$
(the closure of the finite convex combinations in $T$) then
$d_X(p,T)\leq\mbox{diam}(T)$. That is:
$$\mbox{diam}(\mbox{hull}\,T)\leq3 \mbox{diam}(T)\,.$$
\end{itemize}
\end{lemma}
\begin{proof}
The first statement is a straightforward consequence of the definition
of weak convergence. Now let $w,v_1,v_2,\dots,v_k\in T$,\
$\lambda_i\geq 0$, $\sum_{i=1}^k \lambda_i=1$.
It is enough to prove that
$$d_X(\sum^k_{i=1}\lambda_i v_i,w)\leq\sum^k_{i=1}\lambda_i d_X(v_i,w)\,.$$
$$d_X (\sum^k_{i=1}\lambda_i v_i,w)=
\sum^\infty_{n=1}\frac{1}{2^n}|((\sum_{i=1}^k\lambda_iv_i)-w)(f_n)|\leq\\$$
$$\leq \sum^\infty_{n=1}\frac{1}{2^n}\sum_{i=1}^k \lambda_i |(v_i-w)(f_n)|=
\sum_{i=1}^k \lambda_i \sum^\infty_{n=1}\frac{1}{2^n}|(v_i-w)(f_n)|=
\sum^k_{i=1}\lambda_i d_X(v_i,w)\,.$$
\qed \end{proof}
\subsection{Borel-graphings}
Let $X$ be a standard Borel-space and $R\subset X^2$ be a Borel-set
which is a symmetric and irreflexive relation. This structure is called
a Borel-graphing and denoted by $\cal{G}=(X,R)$.
If $xRy$ then we say that $x\sim y$ is an edge of $\cal{G}$.
Let $\Gamma$ be a discrete group acting on a space $X$ in a Borel way
and let $S$ be a generating system of $\Gamma$. Let $x$ and $y$ be
in relation $R$ if $x\neq y$ and $sx=y$ for some $s\in S$. Then
$R$ is a Borel-graphing. The connected components of $\cal{G}$ are
countable graphs on the orbit of the $\Gamma$-action.
We shall use the following result of Kechris, Solecki and Todorcevic 
\cite{KST}: Any Borel-graphing with vertex degree bound $d$ has
a Borel-coloring by $d+1$-colors. That is there is a partition
of $X$ into $d+1$ Borel pieces such that if two points are in the
same piece, then they are not adjacent in $\cal{G}$.

\subsection{Limits of graph sequences}~\label{limitsec}
In this subsection we briefly recall the notion of weak graph convergence
from \cite{BS}.
A graph sequence $\G = \{G_n\}_{n=1}^{\infty} \subset \Gd$ is
  \textit{weakly convergent} if $\lim_{n \to \infty} |V(G_n)| = \infty$
  and for every $r$ and every $\alpha \in U^{r,d}$ the limit
  $\lim_{n\to\infty}p_{G_n}(\alpha)$ exists.

\noindent
 Let $\Grd$ denote the set of all countable, connected rooted
graphs
  $G$ for which $\deg(x) \leq d$ for every $x \in V(G)$.
If $G,H\in\Grd$ let $d_g(G,H)=2^{-r}$, where
$r$ is the maximal number such that the $r$-balls around the roots
of $G$ resp. $H$ are rooted isomorphic. The distance $d_g$  makes
$\Grd$ a compact metric space. Given an
  $\alpha \in U^{r,d}$ let $T(\Grd,\alpha) = \{(G,x) \in \Grd : B_r(x)
  \cong\alpha\}$. The sets $T(\Grd,\alpha)$ are closed open sets
that generate the Borel structure of $\Grd$.

\noindent
We can equip $\Grd$ with an equivalence relation $E$: two rooted
graphs $G,H$ are equivalent ($G \sim_E H$) if they are isomorphic as
graphs (but this isomorphism need not respect the root). It is easy
to see that $E$ is a countable Borel equivalence relation. Also,
convergent graph sequences define a {\bf limit measure} $\mu_{\G}$ on $\Grd$,
where $\mu_{\G}(T(\Grd,\alpha))=\lim_{n\to\infty}p_{G_n}(\alpha)$.
Note however that $\mu$ is not necessary an invariant measure on $\Grd$.
That is why we need the notion of $C$-graphs.

\subsection{$C$-graphs and the space $\CGrd$}~\label{bgraphsec}

 In this subsection we extend our definitions for edge-colored graphs.
A $C$-graph is a graph with edges properly colored by the set
$c_1,c_2,\dots,c_{\ad}$. 
 That is each
edge is labeled by an element of $\{c_1,c_2,\dots,c_{\ad}\}$ and incident
edges are labeled differently. The reason why we use $\ad$ colors will be made 
clear in the last section. We shall denote by $\CGd$ the set of finite
connected $C$-graphs (up to colored isomorphisms).

\noindent
Let $V^{r,d}$ be the set of rooted $(r,d)$-balls edge-colored by
$\{c_1,c_2,\dots,c_{\ad}\}$. Note that we consider two such rooted
$C$-graphs isomorphic if there is a rooted isomorphism between them
preserving the colors. Again, for $CG\in\CGd$ and $\beta\in V^{r,d}$
we define the set $T(CG,\beta)\defeq\{x\in V(CG):\, B_r(x)\cong\beta\}$
and let $p_{CG}(\beta)\defeq\frac{|T(CG,\beta)|}{|V(CG)|}\,.$
We define the statistical distance of two $C$-graphs $CG$ and $CH$ as
$$d_s(CG,CH)=\sum^\infty_{i=1}\frac{1}{2^i}
|p_{CG}(\beta_i)-p_{CH}(\beta_i)|,$$
where $(\beta_1,\beta_2,\dots)$ is an enumeration of all the edge-colored
$(r,d)$-balls, for all $r\geq 1$.

\noindent
Let $\CGrd$ be the isomorphism classes of all connected countable rooted
$C$-graphs with vertex degree bound $d$.
For $\beta \in V^{r,d}$ let us define $T(\CGrd,\beta) = \{x \in
  \CGrd : B_r(x) \cong \al\}$. Again, we have a natural metric
on $\CGrd$. If $X,Y\in \CGrd$ then 
$$d_C(X,Y)=2^{-r}\,,$$
where $r$ is the maximal number such that $B^r(x)\cong B^r(y)$,
where $x$ is the root of $X$, $y$ is the root of $Y$.
The subsets $T(\CGrd,\beta): \beta \in V^{r,d},
  r \in \N$ are closed-open sets and generate the Borel-structure
of $\CGrd$. 

\noindent
Let $\{CG_n\}^\infty_{n=1}$ be a sequence of $C$-graphs.
We say that $\{CG_n\}^\infty_{n=1}$ converges if
for any $\beta \in V^{r,d}$, $\lim_{n\to\infty} p_{CG_n}(\beta)=
\mu(T(\CGrd,\beta))$ exists. In this case $\mu$ is a Borel-measure
on $\CGrd$. We call $\mu$ the limit measure of $\{CG_n\}^\infty_{n=1}$.
\begin{remark} \label{remark2}
We introduce a metric for $P(\CGrd)$ the following way.
Let $\{\beta_1,\beta_2,\dots\}$ be an enumeration of
the classes $V^{r,d},r\geq 1$. Let
$$d_{\CGrd}(\mu,\nu):=
\sum^\infty_{i=1}\frac{1}{2^i} |\mu(T(\CGrd,\beta_i)-\nu(T(\CGrd,\beta_i)|\,.$$
Clearly $d_{\CGrd}$ metrize the weak-topology of $P(\CGrd)$. Also,
$\mu(T(\CGrd),\beta_i))=\int_{\CGrd} 1_{T(\CGrd,\beta_i)}\,d\mu$,
hence Lemma \ref{l21} applies. Note that there is a natural forgetting
map $\Psi:\CGrd\to\Grd$. Again we can consider an enumeration
$(\alpha_1,\alpha_2,\dots)$ of the classes $U^{r,d}$, $r\geq 1$ as
in the Introduction. 
Let
$$d_{\Grd}(\mu,\nu):=
\sum^\infty_{i=1}\frac{1}{2^i} |\mu(T(\Grd,\alpha_i)-\nu(T(\Grd,\alpha_i)|\,.$$
The distance $d_{\Grd}$ metrize the weak-topology of $P(\Grd)$.
By compactness, for any $\epsilon>0$, there exists $\delta>0$ such that
if $d_{\CGrd}(\mu,\nu)<\delta$ then $d_{\Grd}(\Psi_\star \mu,\Psi_\star \nu)
<\epsilon$.

\noindent
Also note that for any finite $C$-graph $CG$ one can associate
a measure on $\CGrd$ concentrated on finitely many points.
Simply define $\mu_{CG}(T(\CGrd,\beta)):=p_{CG}(\beta)$ for any $\beta\in 
V^{r,d}$. The same way for any finite graph $G$ one can associate
a measure $\mu_{G}$ on $\Grd$. Notice that
$$d_s(G,H)=d_{\Grd}(\mu_{G},\mu_{H})\,.$$
Thus the statistical distance is just the distance of the two
associated measures 
$$d_s(CG,CH)=d_{\CGrd}(\mu_{CG},\mu_{CH})\,.$$
In other words, weak convergence of graphs actually means the weak convergence
of the associated measures. Note that if $G$ is a finite graph and
$G'$ is the graph consisting of $2$ disjoint copies of
$G$ then $\mu_G=\mu_{G'}$.
\end{remark}

The equivalence relation $E$ can be extended to $\CGrd$ in an obvious way.
We shall denote this Borel-equivalence relation by $E_C$.
However in this case, we have a natural continuous group action that
defines the equivalence relation. Namely, let $F^2_{\ad}$ be the the free
product of $\ad$ elements of order $2$, with generators
$\{c_1,c_2,\dots,c_{\ad}\}$. The generators act on $\CGrd$ the obvious way :
If $CG$ is an element of $\CGrd$ with root $x$ and $c_i$ is a generator,
then $c_i(CG)$ has the same graph as $CG$ but with root $y$, where $y$ is
the endpoint of the edge colored by $c_i$ (if there is no such edge, then
$CG$ is fixed by $c_i$). Obviously, this action of $F^2_{\ad}$ is
continuous and the orbits of the action are exactly the equivalence
classes. Since the action is continuous, the space of invariant measures
$\I_C$ is compact. Note that any finite $C$-graph $CG$ defines an invariant
point-measure, hence all the limit measures are invariant measures.
This is the advantage of using $C$-graphs instead of colorless graphs.
\begin{remark}
It is worth to mention that connected finite $C$-graphs always
define ergodic measures
and it is easy to construct sequence of connected finite graphs
converging a non-ergodic limit measure. 
\end{remark}

\noindent 
In Section \ref{color}, we shall prove that if $\{G_n\}$ is
a convergent graph sequence then we have a convergent sequence $\{CG_n\}$
of $\ad$-colorings of $\{G_n\}$.

\section{The Three Lemmas}
The three ingredients of our proof of Theorem \ref{mainthm}
are 
the Stability Lemma, the Decomposition Lemma and the
Homogeneity Lemma.
\subsection{The Stability Lemma}
The {\bf edit distance} $ed(G,H)$ of two graphs $G,H$ on the
  same set of vertices is the minimum number of edges that has to be
  deleted from the graphs to make them identical, divided by the
  number of vertices. The edit distance can also be interpreted for
  $C$-graphs the same way, except here we require that the two graphs
  become labeled-identical after the removal of the edges.

\begin{lemma}[The Stability Lemma]~\label{distancelem} {\ }
\begin{enumerate}
\item
For any $\omega > 0$ there is a $\ep > 0$ such that $ed(G,H) < \ep$
implies $d_s(G,H) < \omega$ for all graphs $G,H$.
\item 
For any $\omega > 0$ there is a $\ep > 0$ such that $ed(CG,CH) < \ep$
implies $d_s(CG,CH) < \omega$ for all $C$-graphs $CG,CH$.
\end{enumerate}
\end{lemma}

\begin{proof} We prove the first part here, the proof of the second
  part is basically identical.  Let us suppose that $ed(G,H) <\ep$. Let us fix
a natural number $r$. 
Those vertices that have different $r$-neighborhoods in $G$ and
$H$ must be ``close'' to a deleted edge in one of the graphs. Thus an
upper estimate for the number of such points is
$2\cdot 2dn\ep\cdot d^{r-1}$. Thus if $\alpha \in U^{r,d}$ then
\[|p_G(\alpha) - p_H(\alpha)| \leq 4d^r\ep.\]
Let us choose $i_0$ so that $1/2^{i_0} < \omega/4$. In our enumeration of
neighborhoods $\alpha_i$ let $r_0$ be the largest occuring radius
among the first $i_0$ elements. Finally let us choose $\ep$ to be
smaller than $\frac{\omega}{8 d^{r_0}}$.
With these choices we have
\begin{multline}
d_s(G,H) =  \sum_{i=1}^{\infty}
  \frac{1}{2^i}|p_G(\alpha_i)-p_H(\alpha_i)| = \\
= \sum_{i=1}^{i_0}
  \frac{1}{2^i}|p_G(\alpha_i)-p_H(\alpha_i)| + \sum_{i=i_0+1}^{\infty}
  \frac{1}{2^i}|p_G(\alpha_i)-p_H(\alpha_i)| \leq \\
\leq \sum_{i=1}^{i_0}
  \frac{1}{2^i}4d^{r_0}\ep +
  \sum_{i=i_0+1}^{\infty}\frac{2}{2^i} \leq \omega/2 + \omega/2 = \omega
\end{multline}
 that completes the proof.
\end{proof} \qed

\subsection{The Decomposition Lemma}~\label{splittingsec}

\begin{defin} A \textit{$K$-splitting} of a convergent graph sequence
  $\G = \Gseq$ with limit measure $\mu$ is a collection of $K$ graph
  sequences $\G^i = \Gseq[G^i]\; 1\leq i \leq K$ obtained by removing
  some edges from $\G$ that satisfies the following properties:
\begin{enumerate}
\item $G^i_n$ $(1\leq i\leq K)$ are vertex-disjoint spanned subgraphs of $G_n$.
\item $\lim_{n\to \infty} \frac{|E(G^i_n,G^j_n)|}{|V(G_n)|} = 0$ for
  every $i \neq j$, that is the ratio of removed edges tends to 0. 
\item For any $1\leq i \leq K$,
$\lim_{n\to\infty} \frac {|V(G^i_n)|}{|V(G_n)|} = a_i$ exists.
\item Either $\lim_{n\to\infty} \frac {|V(G^i_n)|}{|V(G_n)|}=0$
or $\{G^n_i\}$ is a convergent graph sequence with limit measure $\mu_i$.
\end{enumerate}
The exact same notion can be defined for sequences of $C$-graphs.
\end{defin}

\begin{propo}~\label{splitting} In a $K$-splitting $\sum_{i: a_i \neq 0} a_i
  \cdot \mu_i = \mu$.
\end{propo}

\begin{proof} Let $H_n = \cup_{i=1}^K G_i$ be the graph obtained from
  $G_n$ by removing the neccessary edges. From part b) of the
  definition of a splitting we have that $ed(G_n,H_n) \to 0$, so by
  the Stability Lemma $d_s(G_n,H_n) \to 0$ hence $H_n$ also converges
  to $\mu$. Now for a fixed neighborhood type $\alpha$ we have
  $p_{H_n}(\alpha) = \sum \frac{|V(G_n^i)|}{|V(H_n)|} \cdot
  p_{G_n^i}(\alpha)$. The left hand side converges to
  $\mu(T(\Grd,\alpha))$, while each term on the right hand sides
  converges to $a_i \cdot \mu_i(\Grd,\alpha)$, even if $a_i=0$. So
  $\mu = \sum a_i \cdot \mu_i$. \qed
\end{proof}

\smallskip

Now let us consider the convex compact space $\I_C$ of invariant
measures on $\CGrd$ and the set of its extremal points, that is the
set of ergodic invariant measures $\EI_C$ (note that $\EI_C$ is non-compact).
\begin{lemma}[The Decomposition Lemma]~\label{partitionlem} 
Let $\CG$ be a  $C$-graph sequence
that converges weakly to an invariant measure $\mu$ on $\CGrd$.  Let
$Z_1,\dots,Z_L$ be a Borel-partition of $\EI_C$.  Then we have a
$K$-splitting $(\CG^1,\dots,\CG^K)$ of $\CG$ such that $\mu_i \in
hull(Z_i)$ whenever $a_i \neq 0$.
 \end{lemma}

\begin{proof}
 Let $Y_i =
  \pi^{-1}(Z_i)$, where $\pi$ is the Borel-surjection
in the Farrell-Varadarajan Theorem.
The sets $Y_i$ are  $E_C$-invariant Borel-subsets of $\CGrd$.
Let $\mu$ denote the limit measure of $\CG$. For any 
natural number $M>0$ one
can approximate the partition $\{Y_i\}^K_{i=1}$ with
a slightly perturbed partition $\{Y^M_i\}^K_{i=1}$
where
\begin{itemize}
\item $\lim_{i\to\infty}
\mu(Y_i\triangle Y^M_i)=0$ for any $1\leq i \leq K$.
\item Each $Y^M_i$ is a closed-open set in the form
$$Y_i^M=\cup^{L_{i,M}}_{j=1} T(\CGrd,\alpha^i_j)\,,$$
where $\alpha^i_j\in V^{M,d}$.
\end{itemize}
That is if $x\in V(CG_n)$ then by looking at the $M$-neighbourhood of $x$
and the colors of its edges one can decide for which
$1\leq i \leq K$; $x\in T(CG_n,Y^M_i)$. Note that we use the notation
$T(CG_n,Y^M_i)$ for $\cup^{L_{i,M}}_{j=1} T(CG_n,\alpha^i_j)$.
Let $\partial Y^M_i$ denote the set of graphs $X$ in $Y^M_i$
such that if we move the root of $X$ to one of its neighbour the resulting
rooted $C$-graph is not in $Y^M_i$. Since $Y_i$ is $E_C$-invariant
set and $\mu(Y_i\triangle Y^M_i)\to 0$, we have $\lim_{M\to\infty}
\mu(\partial Y^M_i)=0$. Note that $\partial Y^M_i$ is
still a closed open set in the form 
$\cup^{R_{i,M}}_{j=1} T(\CGrd,\beta^i_j)$, where 
 $\beta^i_j\in V^{M+1,d}$
Hence $\lim_{n\to\infty} p_{CG_n}(\partial Y_i^M)=\mu(\partial Y_i^M)$.
Observe that $T(CG_n,\partial Y^M_i)$ is the set of vertices $x\in V(CG_n)$
such that $x\in T(CG_n,Y^M_i)$ but $y\notin T(CG_n,Y^M_i)$ for
a neighbour of $x$.
Let
$$Q_M=\{\alpha\in V^{M,d}\,\mid\, \mu(T(\CGrd,\alpha)\neq 0\}$$
$$R_M=\{\alpha\in V^{M,,d}\,\mid\, \mu(T(\CGrd,\alpha)=0\}$$
Since $\CG$ is a convergent $C$-graph sequence
there exists a natural number $n_M$ such that
if $n\geq n_M$ then
\begin{itemize}
\item $p_{CG_n}(\partial Y^M_i)\leq \mu(\partial Y^M_i)+\frac{1}{M}$
 for any $1\leq i \leq K$.
\item $(1-\frac{1}{M})\mu(T(\CGrd,\alpha))<p_{CG_n}(\alpha)<
(1+\frac{1}{M})\mu(T(\CGrd,\alpha))$ for any $\alpha\in Q_M$.
\item $p_{CG_n}(\alpha)\leq \frac{1}{M}\frac{1}{|V^{M,d}|}$, for any
$\alpha\in R_M$.
\end{itemize}
If $n_M\leq n <n_{M+1}$ let us partition $V(CG_n)$ as
$$V(CG_n)=\cup^K_{i=1} V(CG_n^{i,M})=\cup^K_{i=1} T(CG_n,Y^M_i)\,.$$
According to our conditions on $n_M$, if $n>n_M$ then
\begin{equation} \label{becs1a}
\left| \frac{|T(CG_n,\alpha)|}{|V(CG_n)|}-\mu(T(\CGrd,\alpha))\right|\leq 
\frac {1} {M} \mu(T(\CGrd,\alpha))
\quad\mbox{if $\alpha\in Q_M$.}
\end{equation}
\begin{equation} \label{becs1b}
\left| \frac{|T(CG_n,\alpha)|}{|V(CG_n)|}\right|\leq 
\frac {1} {M} \frac{1} {|V^{M,d}|}\quad\mbox{if $\alpha\in R_M$.}
\end{equation}
\begin{equation} \label{becs2}
\left| \sum_{1\leq i,j \leq K} \frac {|E(CG_n^{i,M}, CG_n^{j,M})|}
{|V(CG_n)|}\right| \leq K d\left(\mu(\partial Y^M_i)+\frac{1}{M}\right)
\end{equation}
\begin{lemma}
For any $r>1$ and $\delta>0$ there exists $M>0$ such that for any 
$\beta\in V^{r,d}$ and $1\leq i \leq K$ if $\mu(Y_i)\neq 0$ then
$$\left|\frac{|T(\cgim,\beta)|}{|V(\cgim)|}-\frac{\mu(T(\CGrd,\beta)\cap
Y_i)}{\mu(Y_i)}\right|<\delta\,.$$
\end{lemma}
\begin{proof}
Let $M>r$. Then for any $1\leq i \leq K$ and for any graph $G$, we
have decompositions
$$T(CG,\beta)=\cup^K_{i=1} \cup_{j=1}^{N_{i,M}} T(CG,\alpha^i_j)\,,$$
where $\alpha^i_j\in V^{M,d}$ and
$\cup_{j=1}^{N_{i,M}} T(CG,\alpha^i_j)=T(CG,Y^M_i)\cap T(CG,\beta)$.
Also
$$T(\CGrd,\beta)=\cup^K_{i=1} \cup_{j=1}^{N_{i,M}} T(\CGrd,\alpha^i_j)\,,$$
and $\cup_{j=1}^{N_{i,M}} T(\CGrd,\alpha^i_j)= Y^M_i\cap T(\CGrd,\beta)$.

\noindent
Since $\mu(T(\CGrd,\beta)\cap Y^M_i)\to \mu(T(\CGrd,\beta)\cap Y_i)$
by the estimates (\ref{becs1a}) and (\ref{becs1b}) one can immediately
see that for any $\kappa>0$ there exists $M>0$ such that if $n>n_M$
and $\beta\in V^{r,d}$ then for any $1\leq i \leq K$, where $\mu(Y_i)\neq
0$:
\begin{equation}
\label{becs4}
\left|\frac{|V(\cgim)\cap T(CG_n,\beta)|}{|V(CG_n)|\mu(Y_i)}-
\frac{\mu(T(\CGrd,\beta)\cap Y_i)}{\mu(Y_i)}\right|<\kappa
\end{equation}
That is, if $M$ is large enough then
 if $n>n_M$
and $\beta\in V^{r,d}$ then for any $1\leq i \leq K$, where $\mu(Y_i)\neq
0$:
\begin{equation}
\label{becs5}
\left|\frac{|V(\cgim)\cap T(CG_n,\beta)|}{|V(CG_n)|\frac{|V(\cgim)|}
{|V(CG_n)|}}-
\frac{\mu(T(\CGrd,\beta)\cap Y_i)}{\mu(Y_i)}\right|<\frac{\delta}{3}\,.
\end{equation}

Hence it is enough to show that if $M$ is large enough then
for any $1\leq i \leq K$ and $\beta\in V^{r,d}$
\begin{equation}\label{becsles1}
\left|\frac{(V(CG^{i,M}_n)\cap T(CG_n,\beta))\triangle
T(CG^{i,M}_n,\alpha))}{|V(G_n)|}\right|<\frac{\delta}{3}\,.
\end{equation}
Observe that if $x\in (V(CG^{i,M}_n)\cap T(CG_n,\beta))\triangle
T(CG^{i,M}_n,\beta))$, then 
\begin{equation} \label{becsles2}
B_r(x)\cap T(CG_n,\partial Y_i^M)\neq \emptyset
\end{equation}
Since $\frac{|T(CG_n,\partial Y_i^M)|}{|V(G_n)|}\leq \mu(\partial
Y^M_i)+\frac{1}{M}$, 
for any $\eta>0$ one can choose large enough $M$ so that
$$\frac {|B_r(T(CG_n,\partial Y_i^M))|}{|V(G_n)|}\leq\eta\,.$$
Therefore our lemma follows. \qed \end{proof}

\noindent
By our lemma if $CG^i_n=CG^{i,M}_n$, where $n_M< n \leq n_{M+1}$,
$ \lim_{n\to\infty} \frac { |V(CG^i_n)|}{|V(CG_n)|}\neq 0$ then
$\{CG_n^i\}_{n=1}^\infty$ is convergent and the limit measure
is $\mu_i$, where $\mu_i(U)=\frac{\mu(U\cap Y_i)}{\mu(Y_i)}$.
By the Farrell-Varadarajan Theorem
$$\mu_i=\frac{\int_{Y_i}\pi(x) d\mu(x)}{\mu(Y_i)}$$
that is 
$$\mu_i=\int_{Z_i} p \,d\nu_i(p)\,,$$
where $\nu_i$ is a Borel-probability measure on $Z_i$. Thus $\mu_i$
is the barycenter of $\nu_i$, hence $\mu_i\in hull(Z_i)$. 
\qed \end{proof}

\noindent
\subsection{The Homogeneity Lemma}
\begin{lemma}[The Homogeneity Lemma]
Let $0<\delta<1$, $0<\lambda<1$ be real numbers and let $\delta'<\delta$ be
the constant in Remark \ref{remark2}, that is if for two $C$-graphs
$d_s(CG,CH)<\delta'$ then $d_s(G,H)<\delta$ holds for the underlying graphs.
Let $Z\subset \EI_C$, such that $\diam(Z)<\frac{1}{3}\delta'$.
Let $\{CG_n\}^\infty_{n=1}$ be a sequence of $C$-graphs converging to
a measure $\mu$ supported on $\pi^{-1}(Z)$. Then there exists
$0<\epsilon<\delta$ and $N>0$ such that if $n\geq N$ then
the graph $G_n$ is $(\epsilon,\lambda,\delta)$-homogeneous.
\end{lemma}
\proof
Suppose that the Lemma does not hold. Then we have a sequence of graphs
$\{CG_n\}^\infty_{n=1}$ such that $G_n$ are not
$(\frac{1}{n},\lambda,\delta)$-homogeneous. Thus
$CG_n$ are not $(1/n,\lambda,\delta')$-homogeneous by the choice of $\delta'$. 
Therefore for each $CG_n$ there is a spanned
subgraph $CH_n \subset CG_n$ such that it is not small: $|V(CH_n)|
\geq \lambda |V(CG_n)|$, it has only few edges going out of it:
$|E(CH_n, CG_n\setminus CH_n)| \leq |V(CG_n)|/n$, but still it
is not similar to the big graph: $d_s(CH_n,CG_n) > \delta'$.

We may again choose a subsequence so that $CH_n$, $CG_n
\setminus CH_n$ convergent $C$-graph sequences and 
$\frac{|V(CH_n)|}{|V(CG_n)|}$ is
a convergent sequence of real numbers. The $C$-graphs converge
to invariant measures $\mu_1',\mu_1'' \in \I_C$ while
$\frac{|V(CH_n)|}{|V(CG_n)|}$ converges to some real number $\lambda\leq a \leq
1$. Therefore
we obtained a $2$-splitting of the $CG_n$ sequence. By
Proposition~\ref{splitting} we know that $\mu_1 = a \cdot \mu_1' +
(1-a)\cdot \mu_1''$. The splitting was chosen such that
$d_s(CH_n,CG_n) > \delta'$, so in the limit we have
$d_{\CGrd}(\mu_1,\mu_1') \geq \delta'$. On the other hand since $\mu_1$ was
entirely supported on $\pi^{-1}(Z)$, the same must hold for
$\mu_1'$ and $\mu_1''$.
Since $Z$ has diameter strictly less
than $\frac{1}{3}\delta'$ by Lemma \ref{l21} its convex hull has diameter
less than $\delta'$. This leads to a contradiction. \qed

\section{The proof of Theorem \ref{mainthm}}~\label{mainsec}

Again, let  $\delta' > 0$ be small enough such that for any two $C$-graphs and
their underlying regular graphs $d_s(CG,CH) < \delta'$ implies
$d_s(G,H) < \delta$.
The space $\I_C$ is compact, hence it can be split into $K = K(\delta)$
Borel pieces each of diameter strictly less than $\delta'/3$. Let us
denote by $Z_i$ the intersection of the $i$-th piece with $\EI_C$.  Let
us suppose that for this choice of $K$ there is no good choice for $N$
and $\ep$. That means we have graphs $G_n$ with $|V(G_n)| \geq n$ such
that $G_n$ does not have the desired decomposition into pieces that
are $(1/n, \lambda, \delta)$-homogeneous.

 Let $\{CG_n\}^\infty_{n=1}$ be $\ad$-colorings of $G_n$. By Theorem
 \ref{theorem2} we can assume
that $\{CG_n\}^\infty_{n=1}$
converge to an invariant measure $\mu\in\I_C$. By the Decomposition
Lemma we have a $K$-splitting $\{CG^i_n\}^K_{i=1}$ such that
the limit measure $\mu_i$ of $CG^i_n$ (whenever $a_i$ in the Decomposition
Lemma is non-zero) is entirely supported 
on $\pi^{-1}(\mbox{hull}\,Z_i)\subset \CGrd$.

Since in a $K$-splitting the ratio of deleted edges tends to 0, all
but a finite number of the $CG_n$'s are split by removing less
than $\frac{\delta}{10} |V(CG_n)|$ edges. Let us further remove all edges
from those parts $CG_n^i$ for which
\[\frac{|V(CG_n^i)|} {|V(CG_n)|}\leq \frac{\delta}{10\,d\,K},\]
and put these parts into $CG_n^0$. The number of edges removed in this
step is at most $\frac{\delta}{10\,d\,K} |V(CG_n|$, so in total we still
have not removed more than $\delta |V(CG_n|$ edges. Lastly, the
empty parts together clearly contain less than $\delta |V(CG_n|$
vertices.

We can use this partition of the sequence $\{CG_n\}^\infty_{n=1}$ on the
original graph sequence $\{G_n\}^\infty_{n=1}$  to obtain a candidate for
 the partition
required in the theorem: we didn't remove too many edges, all parts
are big or empty and the empty parts are small altogether. Still our
graphs are counterexamples to the theorem, so for each $n$, one of the
non-empty parts must not be $(1/n,\lambda, \delta)$-homogeneous. This is
in contradiction with the Homogeneity Lemma.

\noindent
Now we prove part (b).  If we have two convergent graph sequences
$G_n, H_n$ for which $d_s(G_n,H_n) \to 0$ then by our Theorem \ref{theorem2}
there exist convergent $d+1$-colorings of them $\{CG_n\}^\infty_{n=1}$
and $\{CH_n\}^\infty_{n=1}$ converging to the same limit measure $\mu$.

 Now we can write
$\mu = \sum a_i \mu_i$ and partition $G_n$ and $H_n$ just as in the
proof of part a). For a fixed index $i$ the $CG_n^i$ will converge to
$\mu_i$ if $a_i > \delta/10dK$ and will be empty otherwise. The
same holds for the $CH_n^i$'s. Hence for large $n$ the same parts will
remain non-empty, and they will both converge to $\mu_i$, while the
ratios
$\frac{|V(CG_n^i)|}{|V(CG_n)|},\frac{|V(CH_n^i)|}{|V(CH_n)|}$
will both converge to $a_i$. 

So if for a fixed $\sigma > 0$ the statment were false, we could
choose a pair of graphs $G_n, H_n$ which would provide a
counterexample for $\tau = 1/n$. However a convergent subsequence of
these graphs would contradict the observation of the previous
paragraph. Thus part b) follows.

\qed
\section{The Coloring Theorem} \label{color}
In this section we prove a folklore conjecture: a convergent graph sequence
can be vertex- resp. edge-colored properly to obtain a convergent sequence
of vertex- resp. edge colored graphs. This result is used in the proof of 
Theorem \ref{mainthm}.
\subsection{$B$-graphs and the space $\BGrd$}

  In this subsection we recall the notion of $B$-graphs from \cite{Elek}.
  Let $B = \{0,1\}^{\N}$ be the Bernoulli space of
  0-1-sequences with the standard product measure $\nu$. A rooted
  $B$-graph is a rooted graph $G$ and a function $\tau_G : V(G) \to
  B$. Two rooted $B$-graphs $G$ and $H$ are said to be isomorphic if there
  exists a rooted isomorphism $\psi : V(G) \to V(H)$ such that
  $\tau_H(\psi(x)) = \tau_G(x)$ for every $x \in V(G)$. The
  set of isomorphism classes of all countable rooted $B$-graphs with degree
  bound $d$ is denoted by $\BGrd$.

\noindent 
Let $U^{k,r,d}$ denote the set of isomorphism classes of rooted 
$r$-balls with
degree bound $d$ and vertices labeled with $k$ (0,1)-digits.
 For a $B$-graph $BG$ and a vertex $x\in V(BG)$ by $B_r^k(x)\in U^{k,r,d}$
 we shall denote
the rooted $r$-ball around $x$ with the labels truncated to the
first $k$ digits. 
For any $\al \in U^{k,r,d}$ and a $B$-graph $BG$ we define the set
  $T(BG,\alpha) \defeq \{x \in V(G): B_r^k(x) \cong_B \alpha\}$ and define
  $p_{BG}(\alpha) \defeq \frac{|T(BG,\alpha)|}{|V(G)|}$.
For $\al \in U^{k,r,d}$ let us define $T(\BGrd,\al) = \{x \in
  \BGrd : B_r^k(x) \cong \al\}$. Again, we have a natural metric
on $\BGrd$. If $X,Y\in \BGrd$ then 
$$d_b(X,Y)=2^{-r}\,,$$
where $r$ is the maximal number such that $B^r_r(x)\cong B^r_r(y)$,
where $x$ is the root of $X$, $y$ is the root of $Y$.
The subsets $T(\Grd,\al): \al \in U^{k,r,d},
  k,r \in \N$ are closed-open sets and generate the Borel-structure
of $\BGrd$. 

\noindent
Let $\{BG_n\}^\infty_{n=1}$ be a sequence of $B$-graphs.
We say that $\{BG_n\}^\infty_{n=1}$ converges if
for any $\al \in U^{k,r,d}$, $\lim_{n\to\infty} p_{BG_n}(\alpha)=
\mu(T(\BGrd,\al))$ exists. In this case $\mu$ is a Borel-measure
on $\BGrd$. We call $\mu$ the limit measure of $\{BG_n\}^\infty_{n=1}$.

\noindent
We can consider the Borel-set $\hb$
 of such rooted $B$-graphs where all the vertex labels
are different. Again, we have the
natural equivalence relation on this $B$-graphs defined by a natural
Borel group action (see \cite{Elek}) making $\hb$ a Borel-graphing
(see Section \ref{notationsec}).
The following proposition is the straightforward consequence
of Proposition 2.2 and Corollary 5.1 of \cite{Elek}. See also
\cite{AL} Example 9.9.
\begin{propo} \label{enlemmam}
Let $\{G_n\}^\infty_{n=1}$ be a convergent graph
sequence.
Let $\{BG_n\}^\infty_{n=1}$ be a uniformly random $B$-labelling
of the vertices of $G_n$. Then with probability $1$ 
$\{BG_n\}^\infty_{n=1}$ converges to an $\I_B$-invariant measure $\mu$
supported entirely on $\hb$.
\end{propo}
\noindent
Let $b_1,b_2,b_3,\dots,b_k$ be elements of $B$ such that
$b_i$ is everywhere zero accept in the $i$-th digit.
We call a $B$-graph a $B_k$-graph if all its vertex labels are
from this set and any two adjacent vertex is labelled differently.
Of course, the limit measure of such graphs are concentrated on
those elements of $\hb$ that are vertex labelled by elements
of $b_1,b_2,b_3,\dots,b_k$. We denote this compact subspace by
$\hat{\textbf{BGr}^k_d}$.

\subsection{Borel-colorings}
According to the theorem of Kechris, Solecki and Todorcevic \cite{KST}
we have a Borel-coloring $c$ of
$\hb$ by $d+1$-colors. Note that this gives us a proper
coloring of the vertices of each rooted $B$-graphs with 
$b_1,b_2,b_3,\dots,b_{d+1}$. Thus each vertex has two labels
one from just $B$, the second is from the set
$b_1,b_2,b_3,\dots,b_{d+1}$.
Now let $\mu$ be a  measure on
$\hb$ then we associate a measure $\mu^c_{d+1}$ on 
$\hb$ concentrated on $\hat{\textbf{BGr}^{d+1}_d}$
Indeed, let $\alpha\in U^{d+1,r,d+1}$ a rooted graph labeled only
by the elements of $b_1,b_2,b_3,\dots,b_{d+1}$. Define
$\mu^c_{d+1}(T(\BGrd,\alpha)):=r_{\alpha}$, where
$r_{\alpha}$ is the $\mu$-measure of such rooted graphs $G$ in $\hb$
such that the $d+1$-coloring of the $r$-ball around the root of $G$
is just $\alpha$. 
\begin{propo}
Let $\{G_n\}^\infty_{n=1}\subset\Grd$ be a convergent graph sequence. Let
$\{BG_n\}^\infty_{n=1}$ be a random $B$-labelling of the graphs converging
to a measure $\mu$ supported entirely on $\hb$ as in
Proposition \ref{enlemmam}. Let $\mu^c_{d+1}$ be the associated
measure. Then there exist  proper vertex colorings of the graphs
by $b_1,b_2,b_3,\dots,b_{d+1}$, such that the resulting $B_{d+1}$-graphs
$\{BG_n^{d+1}\}^\infty_{n=1}$ converge to $\mu^c_{d+1}$.
\end{propo}

\proof
We denote by $H^{d+1,r,d}$ the set of rooted balls properly vertex-colored
by the set $\{b_1, b_2,\dots, b_{d+1}\}$.
Clearly, it is enough to prove that for any $r>1$ and $\e>0$ there exist
proper colorings of $\{G_n\}^\infty_{n=1}$ by $\{b_1, b_2,\dots, b_{d+1}\}$
such a way that if $n$ is large enough then
\begin{equation}
\label{celformula}
|p_{G_n}(\gamma)-\mu^c_{d+1}(T(\hb,\gamma))|<\epsilon
\end{equation}
for any $\gamma\in H^{d+1,r,d}$.
Note that $p_{G_n}(\gamma)=\frac{|T(G_n,\gamma)|}{|V(G_n)|}$, where
$T(G_n,\gamma)$ is the set of vertices $x$ for which the colored
$r$-neighborhood of $x$ is rooted isomorphic to $\gamma$.
From now on we fix an $r>1$ and an $\e>0$.

\noindent
For $\gamma\in H^{d+1,r,d}$ let $W_\gamma$ be the Borel-set of points
in $\hb$ such that $B_r(x)$ is isomorphic to $\gamma$ under the $c$-coloring.
Clearly,
$$\hb=\cup_{\gamma\in H^{d+1,r,d}}W_\gamma$$
form a Borel-partition.
We approximate the Borel-sets $W_\gamma$ by closed-open sets the following
way.
For each $M\geq 1$ let $Y^M_\gamma=\cup^{L_{M,\gamma}}_{i=1}
  T(\hb,\alpha^{M,\gamma}_i)$ where
$\alpha_i^{M,\gamma}\in U^{M,M,d}$ and
\begin{itemize}
\item for any $\gamma\in H^{d+1,r,d}$,
$\mu(W_\gamma\triangle Y^M_\gamma)\leq e(M)$
\item
$Y^M_\gamma\triangle Y^M_{\delta}=\emptyset$ if $\gamma\neq\delta$
\item $\lim_{M\to\infty} e(M)=0$.
\end{itemize}
Then $Y^M_\gamma$ are disjoint Borel-sets and
$\mu(\hb\backslash \cup_{\gamma\in H^{d+1,r,d}}
 (Y^M_\gamma\cap W_\gamma))\leq e(M) |H^{d+1,r,d}|\,.$

\noindent
Note that we immediately have a
(not necessarily proper) coloring $c_M:\hb\to\{b_1,b_2,\dots,b_{d+2}\}$.
Simply, let $c_M(x)=b_i$ if $x\in T(\hb,\alpha^{M,\gamma}_j)$ and the color of
the root of $\gamma$ is $b_i$. Let $c_M(x)=b_{d+2}$
if $x\in\hb\backslash \cup_{\gamma\in H^{d+1,r,d}}
 Y^M_\gamma$.
Thus $c_M$ is a step-function approximation of the Borel-function $c$
depending only on the first $M$ digits of the labels of the $M$-neighborhoods
of the points. Observe that if $x\in \cup_{\gamma\in H^{d+1,r,d}}
 (Y^M_\gamma\cap W_\gamma))$ then $c_M(x)=c(x)$.
Consequently, for any $\gamma\in H^{d+1,r,d}$,
$$\lim_{M \to\infty} \mu^{c_M}_{d+1}(T(\hb,\gamma))=
\mu^{c}_{d+1}(T(\hb,\gamma)).$$

That is if $M$ is large enough, then
$$|\mu^{c_M}_{d+1}(T(\hb,\gamma))-\mu^{c}_{d+1}(T(\hb,\gamma)|<\frac{\e}{10}$$
for any $\gamma\in H^{d+1,r,d}$. It is enough to construct of proper
vertex colorings of the graphs $\{G_n\}^\infty_{n=1}$ such that
if $n$ is large enough then
$$|p_{G_n}(\gamma)-\mu^{c_M}_{d+1}(T(\hb,\gamma))|<\frac{\e}{10}\,,$$
for any $\gamma\in H^{d+1,r,d}$.

\noindent
We call $\alpha\in U^{M,M,d}$ {\bf nice} if $\alpha=\alpha^{M,\gamma}_i$ for
some $\gamma\in H^{d+1,r,d}$. Also, we call $\beta\in U^{M,M+r,d}$ nice
if the following holds:
\begin{itemize}
\item if $s\in B_r(x)$, where $x$ is the root of $\beta$, then $\beta_M(s)$ is
  a nice element of $U^{M,M,d}$.
\end{itemize}
Here, $\beta_M$ denotes the restriction of $\beta$ onto the $M$-ball around
$s$. Notice that if $x\in\hb$ and $x\in T(\hb,\beta)$, where $\beta\in
U^{M,M+r,d}$ is nice, then the $c_M$-coloring of the $r$-neighborhood of $x$
depends only on $\beta$. Let $N^{M,M+r,d}$ be the set of nice balls and
$N_\gamma^{M,M+r,d}$ is the set of those nice balls that the corresponding
$c_M$-coloring of the $r$-neighboorhood of the root is isomorphic to $\gamma$.
Now pick a large enough $M$ such that
\begin{itemize}
\item
$$|\mu^{c_M}_{d+1}(T(\hb,\gamma))-
\mu^{c}_{d+1}(T(\hb,\gamma))|<\frac{\e}{10}\,\mbox{ for any 
$\gamma\in H^{d+1,r,d}$}$$
\item
$$\sum_{\beta\in N^{M,M+r,d}} \mu(T(\hb,\beta))\geq 1-\frac{\e}{10}\,.$$
\end{itemize}
Then for any $\gamma\in H^{d+1,r,d}$
$$0 < \mu^{c_M}_{d+1}(T(\hb,\gamma))-\sum_{\beta\in N_\gamma^{M,M+r,d}}
\mu(T(\hb,\beta))<\frac{\epsilon} {10}\,.$$

Now we construct our proper vertex-colorings of the graph sequence
$\{G_n\}^\infty_{n=1}$.
Let $\{BG_n\}^\infty_{n=1}$ is a sequence of $B$-colorings such that
$\{\mu_{BG_n}\}^\infty_{n=1}$ weakly converges to $\mu$.
Let $M$ be the number determined in the previous paragraph. If
$p\in T(BG_n,\alpha),\alpha\in N_\gamma^{M,M,d}$ for a certain $\gamma$
 then let us color $p$ by the
color of the root of $\gamma$. We color the remaining vertices arbitrarily to
obtain a proper coloring of the underlying graph $G_n$.
Then if $n$ is large enough
\begin{itemize}
\item $|p_{BG_n}(\beta)-\mu(T(\hb,\beta))|<\frac{\e}{10}$ for any $\beta\in
  N^{M,M+r,d}$.
\item The ratio of vertices $x$ for which the $M$-neighbourhood of $x$ is not
  nice is less than $\frac{\e}{5}$.
\end{itemize}
Observe that if $x\in T(BG_n,\beta)$, $\beta\in N_\gamma^{M,M+r,d}$ then in
the vertex colored graph $B_r(x)$ will be isomorphic to $\gamma$. Therefore 
(\ref{celformula}) holds if $n$ is large enough. \qed
\subsection{Edge-colorings}
The goal of this subsection is to prove the following theorem.
\begin{theorem} \label{theorem2}
Let $\{G_n\}^\infty_{n=1}\subset \Gd $
be a convergent graph sequence. Then we have colorings of
the graphs by ${d^2+1 \choose 2}$ colors such that the resulting
graph sequence $\{CG_n\}^\infty_{n=1}$ is convergent as well.
In particular, if $\{G_n\}^\infty_{n=1}$ and $\{H_n\}^\infty_{n=1}$ are
two convergent graph sequences that have the same limit measure, then
there exist colorings of the graph sequences $\{CG_n\}^\infty_{n=1}$ 
and $\{CH_n\}^\infty_{n=1}$ converging to the same measure on $\CGrd$.
\end{theorem}
\proof Let the graph $H_n$ be defined the following way.
\begin{itemize}
\item $V(H_n)=V(G_n)\,.$
\item $(x,y)\in E(H_n)$ if $x\neq y$ and $d_{G_n}(x,y)\leq 2$.
\end{itemize}
Clearly, $\{H_n\}^\infty_{n=1}\subset {\bf Graph}_{d^2}$ is a convergent
graph sequence. Thus, by our previous proposition we have convergent
vertex-colorings of $H_n$ by the colors $\{b_1,b_2,\dots, b_{d^2+1}\}$.
Let $a_1,a_2,\dots,a_{{d^2+1 \choose 2}}$ be the set of pairs 
of different elements
of $\{b_1,b_2,\dots, b_{d^2+1}\}$. Color an edge $(x,y)$ by $a_i$ if
$a_i$ is the pair of colors of $x$ and $y$ in the vertex colorings of $H_n$.
Obviously, we obtain a properly edge-colored graph sequence 
$\{CG_n\}^\infty_{n=1}$ and this sequence is convergent as well. \qed

\vskip 0.2in
\noindent
{\bf Note:} Balazs Szegedy and Omer Angel informed
us that they also proved part (a)
of Theorem 1., using a different argument \cite{Bal}.

G\'abor Elek,

\noindent
Alfred Renyi Institute of the Hungarian Academy of Sciences,
\noindent
POB 127, H-1364, Budapest, Hungary,  elek@renyi.hu

\vskip 0.2in

\noindent
G\'abor Lippner,

\noindent
E\"otv\"os University, Department of Computer Science
\noindent
P\'azm\'any P\'eter s\'et\'any 1/C
H-1117 Budapest
Hungary, lipi@cs.elte.hu

\end{document}